\documentclass[11pt]{amsart}
\usepackage[latin9]{inputenc}
\usepackage{amsmath,amssymb,amscd,amsfonts,amsthm,verbatim}
\usepackage[all,cmtip]{xy}
\usepackage{paralist}
\usepackage[mathscr]{eucal}
\usepackage{color}
\usepackage{graphicx}
\usepackage{pdfsync}
\usepackage[]{fontenc}
\usepackage{enumerate}
\usepackage{pgfplots}
\usepackage{tikz-cd} 
\usepackage[all,cmtip]{xy}

\xyoption{all} \numberwithin{equation}{section}
\newtheorem{thm}[equation]{Theorem}
\newtheorem{lemma}[equation]{Lemma}
\newtheorem{prop}[equation]{Proposition}
\newtheorem{cor}[equation]{Corollary}
\newtheorem*{cor*}{Corollary}

\newtheorem{nota}[equation]{Notation}
\newtheorem*{thm*}{Theorem}
\theoremstyle{definition}
\newtheorem{defn}[equation]{Definition}
\newtheorem{rem}[equation]{Remark}
\newtheorem{ex}[equation]{Example}

\newcommand{\C}{\mathbb C}
\newcommand{\Z}{\mathbb Z}

\DeclareMathOperator{\rk}{rk}

\DeclareMathOperator{\codim}{codim}

\DeclareMathOperator{\sym}{Sym}

\DeclareMathOperator{\cliff}{Cliff}
\DeclareMathOperator{\gon}{gon}
\DeclareMathOperator{\Span}{span}
\DeclareMathOperator{\Ann}{\mathbb{A}nn}

\numberwithin{equation}{section}

\title{The slope of fibred surfaces: unitary rank and  Clifford index}
\author{Enea Riva, Lidia Stoppino}
\thanks{The authors were partially supported MIUR PRIN 2017 ?Moduli spaces and Lie Theory? , by MIUR, Programma
Dipartimenti di Eccellenza (2018-2022) - Dipartimento di Matematica ?F. Casorati?, Universit\`a degli Studi di Pavia and
by INdAM (GNSAGA)}

\begin{document}
\maketitle
\begin{abstract}
We prove new slope inequalities for relatively minimal fibred surfaces, showing an influence of the relative irregularity $q_f$, of  the unitary rank $u_f$ and of the Clifford index $c_f$ on the slope. The argument uses Xiao's method and a new Clifford-type inequality for subcanonical systems on non-hyperelliptic curves.\\
\end{abstract}


\section{Introduction}
Let $f\colon S\rightarrow B$ be a relatively minimal fibred surface. Let $K_f=K_S- f^*K_B$ be its relative canonical divisor and $\chi_f=\chi(\mathcal{O}_S)-\chi(\mathcal{O}_F)\chi(\mathcal{O}_B)$ be its relative Euler characteristic (see Section  \ref{subsec pre fibred}). 
A {\em slope inequality} for the fibred surface is  an inequality of the form:
\begin{equation}\label{general-slope}
K_f^2\geq a\chi_f,
\end{equation}
where $a>0$ is a positive rational number depending on the geometry of the fibration.
The first of this kind of results is the celebrated slope inequality proved by Xiao in \cite{Xiao} and by Cornalba and Harris in \cite{C-H} (see also \cite{stoppino}):
\begin{equation}\label{slope}
K_{f}^2\geq 4\frac{g-1}{g}\chi_f,
\end{equation}
where $g=g(F)$, the genus of a general fibre $F$.
A third proof was given later by Moriwaki in \cite{moriwaki}.
So, here $a$ is an increasing function of $g$.

The rank $g$ vector bundle $f_*\omega_f$ over the base curve $B$ is called  the Hodge bundle of the fibred surface. Note that by Leray's spectral sequence $\chi_f$ coincides with the degree of the Hodge bundle $\deg f_*\omega_f$, and the slope inequality (\ref{slope}) can be rephrased as follows:
$$K_{f}^2\geq 4(g-1)\mu(f_*\omega_f)=2\deg(\omega_F)\mu(f_*\omega_f)=2(K_fF)\mu(f_*\omega_f),$$ 
where $\mu(\mathcal{E})$ as usual denotes Mumford's slope of a vector bundle $\mathcal{E}$ (see Section \ref{subsec: xiao}).

After the seminal papers cited above, several results have been obtained by many authors proving an influence of other natural geometric invariants of the fibred surface on this inequalities.
Call $b:=g(B)$ the genus of the base curve. Let us consider in particular the following invariants. (see also Section \ref{subsec pre fibred}):
\begin{itemize}
\item[-] The relative irregularity $q_f=q(S)-b$.
\item[-]  The unitary rank $u_f$ (\cite{SGT}), i.e. the rank of the unitary summand  $\mathcal {U}$ in the second Fujita decomposition of the Hodge bundle.
\item[-]  The gonality of a general fibre $\gon(f)=\gon(F)$.
\item[-]  The Clifford index  of a general fibre $c_f=\cliff(F)$.
\end{itemize}
In general, $q_f\leq g$ and equality holds if and only if  the fibration is trivial. Many results are known about the relations between $q_f$, $u_f$ and $c_f$.
We recall in particular that for non-isotrivial fibred surfaces, we have 
$q_f\leq \frac{g}{2}+1$  if $b=0$ (\cite{xiao-irregular}), but there exist fibred surfaces  such that $q_f> \frac{g}{2}+1$  \cite{pirola}.
Moreover, we have (\cite{Xiao, CLZ})
$u_f \leq \frac{5g+1}{6}$
and  (\cite{GNB, SGT})
\begin{equation}\label{bound}
q_f\leq u_f\leq g-c_f.
\end{equation}
A sharp bound is not known, it is predicted by the modified Xiao conjecture \cite{GNB, SGT}.

The first two invariants satisfy inequality $u_f\geq q_f$ (see Section \ref{subsec: xiao}), but there are fibred surfaces where strict inequality hold: thanks to the results of Catanese and Dettweiler (\cite{CD} and \cite{CD2}) we know that these fibred surfaces are -modulo base change- precisely the ones having unitary summand  $\mathcal U$ with infinite monodromy (see \cite{SGT} and Section \ref{subsec pre fibred}, Remark \ref{rem: uf}). 

The relative irregularity has a clear geometric meaning as the dimension of a fixed abelian variety which is the image  of  the Jacobians of all the smooth fibres via the homomorphism induced by the inclusion \cite{Beau}. On the other hand, the unitary rank has a more elusive meaning (see \cite{SGT}, \cite{TG}).

The question of an increasing bound depending on $q_f$ dates back to the original paper of Xiao, where he proves that, if $q_f>0$, then $K_f^2\geq 4\chi_f$ and, moreover, that if $K_f^2= 4\chi_f$, then $q_f=1$ (\cite[Cor.2, Thm.3]{Xiao}).

Let us consider the gonality and the Clifford index of the fibration: recall that (Section \ref{subsec: canonical})
$$\gon(f)-3\leq c_f\leq \gon(f)-2. $$

It was proved by Konno in \cite{konno-nonhyp} (see also \cite{C-H} and \cite{stoppino}) that if a non-locally trivial fibred surface satisfies equality in (\ref{slope}), then necessarily it is hyperelliptic (i.e. with minimal gonality 2 i.e. Clifford index 0). Thus, one naturally would expect that there exists a function in (\ref{general-slope}) increasing with $\gon(f)$, with some genericity assumption needed, as observed by Barja and the second author in \cite[Remark 3.6]{BS-trigonal}. 

Some results are known for small gonality (\cite{konno-trig, BS-trigonal, BZ, C-S, Eno}).
A stunning approach using relative Koszul sequences is proposed in \cite{konno-cliff}, proving in particular that for odd $g$ and general gonality $\lfloor\frac{g+3}{2}\rfloor=\frac{g+3}{2}$, hence Clifford index  $\lfloor\frac{g-1}{2}\rfloor= \frac{g-1}{2}$, we have:
\begin{equation}\label{general}K^2_{f} \geq 6\frac{g-1}{g+1}\chi_f.\end{equation}
Another step towards an answer to both problems was given by Barja and the second named author in \cite{SB} with the following result:
\begin{equation}\label{eq: proj}K^2_{f} \geq 4\frac{g-1}{g-\lfloor m/2 \rfloor}\chi_f,\end{equation}
where $m=\min\{ q_f, c_{f} \}.$ This bound is interesting if both invariants are big with respect to $g$: this can very well happen, as proved  in loc. cit. by providing several examples.

Very recently, Lu and Zuo introduced a yet another very natural technique, using the relative multiplication map
$$\sym^2f_*\omega_f\longrightarrow f_*\omega_f^{\otimes 2},$$
combined with Xiao's method.
Thanks to this technique, the two authors were able to improve \cite{SB} in both directions.
Firstly, they obtain an inequality with $a=a(g,q_f)$ increasing with the relative irregularity \cite{LZ2}:
\begin{equation}\label{eq: LZ1}K_f^2\geq 4\frac{g-1}{g-q_f/2}\chi_f.\end{equation}
Moreover, they proved in  \cite{LZ3} the following: if a general fibre of $f$ is general in the $k$-gonal locus $\mathcal{D}_k$ in $\mathcal{M}_g$ and $g\geq (k-1)^2$, then 
\begin{equation}\label{eq: LZ2}K_f^2\geq \frac{(5k-6)(g-1)}{(k-1)(g+2)}\chi_f.\end{equation}
No bound is known -to our knowledge- involving the unitary rank $u_f$. 

\bigskip

In this paper we prove new bounds depending increasingly on $c_f$, $q_f$ and $u_f$.
Let us summarize here the main results obtained.
 \begin{thm*}[Theorems \ref{sim} and \ref{Y}]\label{MAIN} Let $f\colon S\rightarrow B$ be a relatively minimal fibred surface of genus $g\geq 2$; let $m:=\min\{ q_{f},c_{f}\}$. The following inequalities hold:
\begin{equation}\label{eq: m}
K_{f}^{2}\geq 2\frac{2g-2-m}{g-m}\chi_{f},
\end{equation}
\begin{equation}\label{eq: ueffe}
K_{f}^{2} \geq \left\{ \begin{array}{ll}
                                      2\frac{2g-2-u_{f}}{(g-u_{f})}\chi_f \ \ \ \ \ \ \ \ \ \ \ \ \text{ if $\,u_{f}\leq c_{f}$;}  \\
                                      2\frac{(2g-2-c_{f})(g-1-u_{f})}{(g-1-c_{f})(g-u_{f})}\chi_f \ \ \ \text{ if $\,u_{f}\geq c_{f}.$}
                                     \end{array} \right.  
\end{equation}
\end{thm*}
\begin{rem}\label{rem: comparison}
 Let us compare our results with the known results.
 \begin{itemize}
\item The first inequality (\ref{eq: m}) improves (\ref{eq: proj}) and, more importantly, is greater than (\ref{eq: LZ1}) in case $q_f\leq c_f$. 
On the other hand  in case $q_f\geq c_f$, inequality (\ref{eq: m}) gives a bound increasing  with the Clifford index. Inequality (\ref{eq: LZ2}) can be better, but 
inequality (\ref{eq: m}) holds also  when (\ref{eq: LZ2}) is not applicable: no  genericity assumptions is needed, nor assumptions on $g\gg m$. Moreover, for $m$ big, or $u_f$ and $c_f$ close to $\frac{g-1}{2}$, the bound of inequalities (\ref{eq: m}) and (\ref{eq: ueffe}) becomes close to 6 (see Remark \ref{rem: asynt} below).
\item Inequalities (\ref{eq: ueffe}) are the first known slope inequalities showing an influence of $u_f$. 
\item Inequalities (\ref{eq: ueffe}) are of particular interest in view of the fact cited above that $u_f$ can be strictly bigger than $q_f$. In Section \ref{example}, following \cite{CD3}, we give a first example of a fibred surface where the second inequality is new. 
This fibred surface has invariants $g=6$, $q_f=0$, $c_f=u_f=2$, and is not bielliptic. The bound of (\ref{eq: ueffe}) is
$K_S^2\geq 4\chi_f$, while the other previously known bounds are strictly smaller or not applicable. 
\end{itemize}
\end{rem}
\begin{rem}\label{rem: asynt}
Note that  all our bounds are asymptotically close to $4$ for $g\gg 0$, and this is natural in view of all the known examples and conjectures. But when $m$ is big with respect to $g$, the slope gets bigger, going asymptotically to $6$. 
Let us observe that for odd genus, if the Clifford index is maximal $\cliff(f)=\lfloor\frac{g-1}{2}\rfloor$ and if  $q_f\geq\frac{g-1}{2}$ (\ref{eq: m}) becomes Konno's bound (\ref{general}).
For Clifford index  (hence gonality) close to $\frac{g-1}{2}$, yet not maximal, these bounds are new.
\end{rem}
\bigskip

Our arguments make use of Xiao's method (Section \ref{subsec: xiao}). 
Basically, Xiao's technique works as follows: given a subsheaf of the Hodge bundle $\mathcal G\subseteq f_*\omega_f$, consider the linear sub-canonical system $\mathcal{G}\otimes\mathbb C(t)\subseteq H^0(F, K_F)$ induced on a general fibre $F=f^{*}(t)$. If one has a lower estimate on the ratio of degree over projective dimension of the linear subsystems of $\mathcal{G}\otimes\mathbb C(t)$, then the method produces an inequality of the form $K_f^2\geq b \deg(\mathcal{G})$, where $b$ is a positive number depending on the lower estimate above. See Section \ref{subsec: xiao} and Theorem \ref{relative} for precise statements. 
Taking as $\mathcal G$ the whole Hodge bundle, Clifford's Theorem (see Section \ref{subsec: canonical}) gives the slope inequality (\ref{slope}).

It is thus very natural to try and apply Xiao's method to the ample summand $\mathcal A$ of the second Fujita decomposition of the Hodge bundle (\ref{second fujita}), as $deg\mathcal{A}=\chi_f$. In \cite{SB-survey}, \cite{SB} the analog approach is discussed with the positive summand of the first Fujita decomposition (\ref{first fujita}). One of the  difficulties with these approaches is that there seems to be no control on the base locus of the linear sub-canonical systems induced by $\mathcal A$ on the general fibres of $f$, neither on the linear stability (ref. Section \ref{subsec: canonical}) of this system. 

However, one can still look to a lower bound for the ratio of degree over projective dimension of the linear subcanonical systems that improves Clifford's bound $2$.

This is what we do in our paper, obtaining  a new Clifford-type inequality for subcanonical systems over a  non-hyperelliptic curve $C$, only depending on the codimension and on the Clifford index of $C$. This gives also the desired control on the base locus of the subcanonical systems.

\begin{thm*}[Theorem \ref{Lemma1}]
Let $C\subseteq \mathbb{P}^{g-1}$ be a canonical non-hyperelliptic curve. Let $V\subseteq H^{0}(C,\omega_{C})$  be a linear subspace of codimension $k\leq g-2$.
Then for any $W \subseteq V $ subspace of dimension $\dim W \geq 2$, we have:
$$\frac{\deg |W|}{\dim |W|}\geq \frac{2g-2-m}{g-m-1},$$
where $m:=\min \{ k, \cliff(C) \}$.
\end{thm*}
 Although the motivation in this paper is to apply Xiao's technique, we believe that this result is interesting on its own. The arguments are of genuine geometric classical flavour. 
 
 The above result implies a stability result, as follows (see Section \ref{subsec: canonical} for the definitions).
\begin{cor}[Corollary \ref{corollo}]
Given $V\subseteq H^0(C,\omega_{C})$ a vector subspace of codimension $k$ and dimension $\geq 2$, with $k\leq \cliff(C)$. Then
$\deg|V|\geq 2g-2-k$, i.e. the base locus of $|V|$ has degree smaller or equal to $k$.
If $\deg|V| =2g-2-k$, then $|V|$ is linearly semistable and in particular it is Chow semistable.
\end{cor}
\begin{rem}
This result should be compared also to  \cite{MS}, where linear stability of linear systems on curves  is discussed in relation to the Clifford index.
\end{rem}
\begin{rem}
The slope inequalities have applications both to the geography of surfaces of general type (see for instance \cite{pardini}) and to the ample cone of the moduli space of curves (see for instance \cite{moriwaki} and \cite{GKM}). These perspectives were the original point of view of Xiao and of Cornalba and Harris respectively. In the last years, many authors have treated the case of slope inequalities of fibrations over curves with total space of higher dimension (see for instance \cite{barja-ine}). 
\end{rem}

The paper is organized as follows. In Section  \ref{sec: ine}, after some preliminaries on canonical curves and linear stability, we prove the main Clifford-type result for non-complete sub-canonical systems on non-hyperelliptic curves. We then discuss some stability consequence and give some natural examples.

In Section \ref{sec: xiao} we start by reviewing in \ref{subsec pre fibred} some basic results on fibred surfaces and their relative invariants. Then in \ref{subsec: xiao} and \ref{subsec: mainxiao} we give a review of the main theorems of Xiao's technique, in the form needed for our arguments. We state Xiao's method for fibred surfaces in full generality, following Konno's and Barja's papers, for any locally free subsheaf $\mathcal {G}$ of $f_*\mathcal{O}_S(D)$, where $D$ is a nef divisor on $S$. 

The proof of the main inequalities is carried on in Section \ref{sec: main}. 

In the last Section, following Catanese and Dettweiler's examples, we provide a first example of a fibred surface such that the inequality in the main theorem involving $u_f$ is new.

\begin{nota}
\upshape{We work over the complex field $\C$. All varieties, unless otherwise stated, are assumed to be smooth and projective. Given a variety $X$ and a divisor $D$ on $X$, $H^0(X,D)$ means as usual $H^0(X, \mathcal{O}_X(D))$.}
\end{nota}


\subsubsection*{Acknowledgements} Firstly we wish to thank Miguel \`Angel Barja for countless invaluable conversations on the subject, and for a careful revision of a first draft of the paper. We also are grateful to Marc Coppens for pointing out a gross misattribution in a first draft of the paper. Last but not least, we acknowledge the anonymous referee for having considerably improved the paper.

\section{Clifford-type inequalities for sub-canonical systems}\label{sec: ine}

\subsection{Preliminaries on canonical curves}\label{subsec: canonical}

Let  $C$ be a smooth projective curve of genus $g(C)=g\geq 2$, and let  $K_{C}$ (resp. $\omega_{C}$) its canonical divisor (resp. line bundle). Let 
\[
\phi_{K} : C \longrightarrow \mathbb{P}(H^0(C,\omega_{C})^\vee)\cong\mathbb{P}^{g-1} 
\] 
be its canonical morphism. Assume that $C$ is non-hyperelliptic, i.e that $\phi_{K}$ is an embedding.
Often, with abuse of notation, we identify $C$ and its points with the corresponding canonical image.



Given a linear subspace $V\subseteq H^0(C,\omega_{C})$, let us consider:
$$Ann(V):=\{ \theta \in H^0(C,\omega_{C})^{\vee}\ |\ \ \theta(v)=0 \ \ \ \forall v \in V\}\subseteq H^{0}(C,\omega_{C})^{\vee}.$$
We call this subspace \textit{annihilator} of $V$.\ Let $\Ann(V)=\mathbb{P}(Ann(V))\subseteq \mathbb{P}(H^0(C,\omega_{C})^\vee)$ be its projectivisation. Observe that the dimension of $\Ann(V)$ is the codimension of $V$ minus one.
\begin{defn}
Given an effective divisor $D$ on $C$, its projective span is 
$$\Span(D)=\Span\phi_{K}(D):=\Ann(H^0(C, \omega_{C}(-D))\subseteq \mathbb{P}(H^0(C,\omega_{C})^\vee)$$
\end{defn} 
\begin{ex} Given a point $p\in C$, $\Span\phi_{K}(p)=\{p\}$, while $\Span\phi_{K}(2p)$ is the line tangent to $C$ in $\mathbb P^{g-1}$,  $\Span\phi_{K}(3p)$ is the osculating plane to $C$, and so on. For $n$ distinct points $p_1,\ldots , p_n$ on $C$ if we call $D=p_1+\ldots +p_n$, we have that 
$\Span\phi_{K}(D)$ coincides with the linear projective span of the points in $\mathbb{P}^{g-1}$.
\end{ex}
\begin{thm}[Geometric version of Riemann-Roch \cite{ACGH1}]\label{grr} Given an effective divisor $D$  of degree $d$ on a smooth non-hyperelliptic curve $C$ of genus $g\geq 2$, we have:
$$\dim (\Span (D))= \dim (\Span\phi_{K}(D))= d-1-\dim|D|=d-h^0(C,D).$$
\end{thm}

Given a linear subspace $V\subseteq H^0(C,\omega_{C})$ consider the scheme-theoretic intersection 
$D_{V}:=\Ann(V)\cap C. $
The divisor  $D_{V}$  is the base locus of the linear system $|V|$ seen as a subsystem of $|K_C|$ of degree 2g-2.
Observe that  the evaluation map of  $V$ is surjective onto $\omega_{C}(-D_{V})$.

%
\begin{defn}(Gonality)
The \textit{gonality} $\gon(C)$ of $C$ is the following integer:
$$ \gon(C):=\min \{\deg(\pi) | \ \ \pi : C\rightarrow \mathbb{P}^1\ \  \text{is a surjective morphism} \}=\min \{m | \,\, \exists\, \, g^{1}_{m}\mbox{  over }C\}.$$
\end{defn}
\begin{defn}(Clifford index)
Given a curve $C$ of genus $g\geq 4$, we define its Clifford index $\cliff(C)$ as:
$$\cliff(C):=\min\{ \deg(D)-2(\dim|D|) \ \ | \ h^{0}(C,D)\geq 2,\ h^1(C,D)\geq 2 \}. $$
In case $g=2,3$ we define the Clifford index as follows:
\begin{itemize}
\item if $g=2$, $\cliff(C):=0$;
\item if $g=3$, $\cliff(C):=0$ (resp. $1$) if $C$ is hyperelliptic (resp. trigonal).
\end{itemize}
\end{defn}

For every divisor $D$ such that $h^0(C,D)\geq 2$ and $h^1(C,D)\geq 2$, we say that $D$ {\em contributes} to the Clifford index.
\begin{rem}
Clifford's Theorem (\cite[pp.107-108]{ACGH1}) is equivalent to the following statement: for any curve $C$ of genus $g\geq 2$, $\cliff(C)\geq 0$ and equality holds if and only if $C$ is hyperelliptic.
\end{rem}

Gonality and Clifford index are well studied invariants. We briefly recall some classical results about them.
We have the following upper bounds: 
\begin{eqnarray*}
\gon(C) \leq \lfloor \frac{g+3}{2}\rfloor , \ \ \ \
\cliff(C)\leq \lfloor \frac{g-1}{2} \rfloor,
\end{eqnarray*}
with equality holding for a general curve in $\mathcal{M}_{g}$.
Gonality also has a very natural geometric interpretation via Geometric Riemann-Roch Theorem:
\begin{prop}
For every effective divisor $D$ over $C$, 
$$\dim (\Span (D))\leq  \deg(D)-1.$$
If $\dim \Span (D)< \deg(D)-1$, then $\deg D\geq \gon(C)$. If on the other hand $k$ is an integer greater or equal to $\gon (C)$, then there exists a divisor $D$ of degree $\deg D=k$ with $\dim \Span D< \deg D-1$.
\end{prop}
\begin{proof}
The first inequality is straightforward from Geometric Riemann Roch \ref{grr}. Suppose now that $\dim \Span (D)< \deg(D)-1$; by Riemann Roch again $h^0(C, D)\geq 2$. So there exists a linear subspace $V\subseteq H^0(C,D)$ of dimension $2$ and degree $\leq \deg D$. Thus $\gon(C)\leq \deg D$. The other implication is immediate. 
\end{proof}
\begin{rem}
For example, a non-hyperelliptic curve $C$  is  trigonal if and only if there exist three collinear points on $C$, 
a curve $C$ is  $4$-gonal (i.e.  $\gon(C)=4$) if and only if  every three points $p_{1},p_{2},p_{3}$ of $C$ are not collinear, but there exist a 4-uple of points of $C$ that spans a plane.
\end{rem}
\begin{rem}
The following inequalities hold, proved by Coppens and Martens \cite{Martens}:
\begin{equation}\label{sopra}
\gon(C)-3\leq \cliff(C)\leq \gon(C)-2. 
\end{equation}
Moreover,  for a general curve $C$ in the locally closed subset of curves in the moduli space of gonality $\gon(C)$, it holds equality on the right (see \cite{Ballico}).
\end{rem}
Eventually, we recall the following definition due to Mumford \cite{Mumford}.
\begin{defn}(Linear (semi)stability)
A linear system $|V|$ over $C$ is {\em linearly stable} (resp. semistable) if for every linear subsystem  $|W|\subseteq |V|$ we have:
$$ \frac{\deg |W|}{\dim |W|}>\frac{\deg |V|}{\dim |V| }\,\,\, (\text{resp. $\geq$})$$
\end{defn}

\begin{rem}
Let us make some remarks.
\begin{itemize}
\item  The linear system $|V|$ and its linear subsystems $|W|$ are not necessarily complete;
\item If $|V|\subseteq |L|$ has a non zero base locus $D$, then the linear subsystem:
$$V(-D):= V\cap H^0(C,L-D)  $$
destabilizes it, because $\deg|V(-D)|<\deg|V|$ but $\dim|V(-D)|=\dim|V|$. So, systems with base points are linearly unstable.
\item\label{bpf} Clearly the definition could be modified considering only base-point free subsystems of $|V|$. 
\item Again, Clifford's theorem can be rephrased saying that the canonical system on a curve is linearly semistable, and it is stable if and only if the curve is non-hyperelliptic.
\item Linear stability was introduced by Mumford in order to develop a simple method to prove GIT stability results, indeed, it is proven in \cite{Mumford} that linear semistability implies Chow stability and in \cite{ACGH2} that linear stability implies Hilbert stability.
\end{itemize}
\end{rem}

\subsection{The main result}
A linear system $|V|$ is linearly stable if its ratio $\deg|V|/\dim|V|$ bounds from below the ratio $d/r$ for any $g^r_d \subseteq |V|$.\
Changing point of view, given a linear system on a curve, one can ask for a lower bound for this ratio $d/r$ possibly lower than the original ratio $\deg|V|/\dim|V|$.
This is what we do for canonical subsystems of non-hyperelliptic curves, obtaining a bound depending on the codimension and on the Clifford index of the curve.
\begin{thm}\label{Lemma1}
Let $C\subseteq \mathbb{P}^{g-1}$ be a canonical non-hyperelliptic curve. Let  $V\subseteq H^{0}(C,\omega_{C})$ a linear subspace of codimension $k\leq g-2 $.
Then for any $W \subseteq V $ subspace of dimension $\dim W \geq 2$, we have:
$$\frac{\deg |W|}{\dim |W|}\geq \frac{2g-2-m}{g-m-1}. $$
where $m:=\min \{ k, \cliff(C) \}$.
\end{thm}
\begin{proof}
For any $W \subseteq V$ we have the evaluation morphism:
$$ W\otimes \mathcal{O}_{C} \twoheadrightarrow \omega_{C}(-D_{W}), $$
where $D_{W}:=\Ann(W)\cap C$ is the base locus of $|W|$ seen as a subsystem of the canonical series.

We begin by considering the case $m=k$. 
\begin{lemma}\label{Lemmino}
If $k\leq \cliff(C)$, then $\deg|V|\geq\deg(\omega_{C}(-D_{V}))\geq 2g-2-k$, i.e. $\deg D_{V}\leq k$.
\end{lemma}
\begin{proof}
We split the proof of the lemma in two cases:
\begin{itemize}
\item If $h^{0}(C,D_{V})\geq 2$, since $h^{0}(C,\omega_{C}(-D_{V}))\geq \dim V \geq 2$, then both $D_{V}$ and $\omega_{C}(-D_{V})$ contributes to the Clifford index of $C$, so we have: 
$$\deg(\omega_{C}(-D_{V}))\geq 2(h^{0}(C,\omega_{C}(-D_{V}))-1)+\cliff(C)\geq 2(g-k-1)+k =2g-2-k, $$
as wanted.
\item If $h^{0}(C,D_{V})=1$, by the geometric version of Riemann-Roch (Theorem \ref{grr}), we have that:
 $$\dim(\Span (D_{V}))=\deg D_{V} -h^{0}(C,D_{V})=\deg D_{V}-1 .$$
 Now, $\Span(D_{V})\subseteq \Ann V$ by construction, and
 $$\dim \Ann(V)=g-1-\dim(V) = g-1-(g-k)=k-1 .$$
Therefore, we can conclude that $\deg D_{V}\leq k$, and the claim is proven. 
\end{itemize}
\end{proof}
Let's go back to the proof of Theorem \ref{Lemma1}.
Let $W\subsetneq V$, with $\dim W \geq 2$. As done for Lemma \ref{Lemmino}, we analyze the two following cases:
\begin{itemize}
\item[(i)] If $h^{0}(C,D_{W})\geq 2$, hence $D_{W}$ contributes to $\cliff(C)$ since $h^1(C,D_{W})=h^0(C,\omega_{C}(-D_{W}))\geq \dim W \geq 2$, then: 
$$\deg \omega_{C}(-D_{W})\geq 2(h^{0}(C,\omega_{C}(-D_{W}))-1)+\cliff(C)\geq 2(\dim W -1) +k .$$
Hence:
$$\frac{\deg |W|}{\dim |W|}= \frac{\deg \omega_{C}(-D_{W})}{\dim |W|}\geq 2+\frac{k}{\dim |W|} \geq 2+\frac{k}{\dim |V|}= \frac{2g-2-k}{g-k-1}, $$
as wanted.
\item[(ii)] If $h^{0}(C,D_{W})=1$ we can conclude $\deg D_{W}\leq \dim (\Ann W)+1$ as in the proof of lemma \ref{Lemmino}. Setting $k_{W}:=\dim(\Ann W)+1=\codim(W) $, we have:
$$\frac{\deg|W|}{\dim|W|}\geq \frac{2g-2-\deg D_{W}}{g-k_{W}-1}\geq \frac{2g-2-k_{W}}{g-k_{W}-1}. $$

Since $W\subseteq V$ we can conclude that $k_{W}\geq k$.\\
Now, consider the function:
\begin{equation}\label{funzione}
 f\colon[0,g-1]\rightarrow \mathbb{R} \ \ \ f(t):=\frac{2g-2-t}{g-t-1}.
 \end{equation}
As
$$f'(t)=\frac{g-1}{(g-t-1)^{2}}>0 \ \ \ \ \forall t \in [0,g-1],$$
we have that $f$ is monotonically strictly increasing. So, since $k_{W}\geq k$, we obtain:
$$\frac{\deg|W|}{\dim|W|}\geq f(k_{W})\geq f(k)=\frac{2g-2-k}{g-k-1},$$ 
as wanted.
\end{itemize}

Let us now treat the case $k\geq \cliff(C)=:c$. We prove that for any $W \subseteq V$, with $\dim W\geq 2$:
$$\frac{\deg |W|}{\dim |W|}\geq \frac{2g-2-c}{g-c-1}. $$

Like we did above, we focus on two cases:
\begin{itemize}
\item[(i)] if $h^0(C,D_{W})\geq 2$, then $D_{W}$ contributes to the Clifford index since 
$$h^1(C,D_{W})=h^0(C,\omega_{C}(-D_{W}))\geq \dim W \geq 2.$$
So we have that
$$\deg(\omega_{C}(-D_{W}))\geq 2(h^0(C,\omega_{C}(-D_{W}))-1)+c \geq 2\dim |W|+c . $$
Then it follows that:
$$\frac{\deg|W|}{\dim |W|}\geq 2 +\frac{c}{\dim |W|}\geq 2+ \frac{c}{g-c-1}=\frac{2g-2-c}{g-c-1}.$$
\item[(ii)] If otherwise $h^0(C, D_{W})=1$, then as in the previous case we can conclude:
$$\deg D_{W}\leq k_{W} $$
and since $k_{W}\geq k\geq c$, exploiting the monotonicity of the function $f$:
$$\frac{\deg|W|}{\dim |W|}\geq\frac{2g-2-\deg D_{W}}{g-1-k_{W}}\geq \frac{2g-2-k_{W}}{g-1-k_{W}}= f(k_{W})\geq f(c)=\frac{2g-2-c}{g-c-1}. $$ 
\end{itemize}
\end{proof}
\begin{rem}
Theorem \ref{Lemma1} above is not a linear stability result for the system $|V|$ unless $k\leq\cliff(C)$ and $\deg |V| = 2g-2-k $, i.e. $D_V$ is of maximal degree according to Lemma \ref{Lemmino}. 
\end{rem}
\begin{cor}\label{corollo}
Let $V\subseteq H^0(C,\omega_{C})$ be a vector subspace of codimension $k$, with $k\leq \cliff(C)$. If
$$\deg|V|=2g-2-k$$
then $|V|$ is linearly semistable. In particular the morphism induced on $C$ is Chow semistable.
\end{cor}
\begin{proof}
Let  $W\subseteq V$. Let $h\geq k$ be the codimension of $W$ in $H^0(C,\omega_{C})$.
By Lemma \ref{Lemma1} we have that, for $\overline m=\min\{\cliff(C), h\}$, 
$$\frac{\deg |W|}{\dim |W|}\geq \frac{2g-2-\overline m}{g-h-\overline m}.$$
Now, $\overline m\geq k$, and we are done by the monotonicity of the function ${f}$ defined in (\ref{funzione}):
$$\frac{\deg |W|}{\dim |W|}\geq \frac{2g-2-\overline m}{g-\overline m-1}=f(\overline m)\geq f(k)= \frac{2g-2-k}{g-k-1}= \frac{\deg |V|}{\dim |V|}. $$
\end{proof}

\begin{ex}
Given $k\leq \cliff C$ points $p_{1},\ldots,p_{k}$ on $C$ in general position, clearly the system $|\omega_{C}(-p_{1}\ldots-p_{k})|$ satisfies the assumptions of Corollary \ref{corollo}, as 
$$\deg(\omega_{C}(-p_{1}\ldots-p_{k}))=2g-2-k \,\, \mbox{ and } \,\,\,h^0(C, \omega_{C}(-p_{1}\ldots-p_{k}))=g-k.$$
\end{ex}
\begin{ex}
We see here that indeed for {\em any} set of $k\leq \cliff(C)$ points on $C$, the system $|\omega_{C}(-p_{1}\ldots-p_{k})|$ satisfies the assumptions of Corollary \ref{corollo}.
Indeed, we claim that 
$$h^0(C, p_{1}+\ldots+p_{k})=1. $$
Assume by contradiction that $h^0(C, p_{1}+\ldots+p_{k})\geq 2$: we would have a $g^1_d$ on $C$ with $d\leq k$ hence 
$$\gon(C)\leq d, $$ 
but from the above mentioned result (\ref{sopra}) we obtain:
$$k+2\leq \gon(C)\leq d\leq k,$$
which gives a contradiction.
From Riemann-Roch theorem 
$$h^0(C, \omega_{C}(-p_{1}\ldots-p_{k}))= 2g-2-k +1-g +h^0(C, p_{1}+\ldots+p_{k})=g-k.$$
Hence the linear series $|\omega_{C}(-p_{1}\ldots-p_{k})|$ satisfies the hypothesis of Corollary \ref{corollo}, so it is linearly semistable.  
\end{ex}

\section{Xiao's method for subsheaves}\label{sec: xiao}

\subsection{Preliminaries on fibred surfaces}\label{subsec pre fibred}

\begin{defn}
We call \textit{fibred surface} or sometimes simply {\em fibration} the data of a morphism $f\colon S\rightarrow B$ from a smooth projective surface $S$ to a smooth projective curve $B$ which is surjective with connected fibres.
\end{defn}
We denote with $b=g(B)$  the genus of the base curve. A general fibre $F$ is a smooth curve and its genus $g=g(F)$ is by definition the genus of the fibration. From now on, we consider fibrations of genus $g\geq 2$.

Let $K_{f}:=K_{S}-f^{*}K_{B}$ (resp. $\omega_{f}:=\omega_{S}\otimes(f^{*}\omega_{B})^\vee$) the relative canonical divisor (resp. line bundle).
Recall that given a surface $S$ a $(-1)$-curve is a non-singular rational curve $C\subseteq S$ such that $C^2=-1$.
We say that $f$ is \textit{relatively minimal} if it does not contain any $(-1)$-curves in its fibres. This condition is equivalent to $K_{f}$ being a relatively nef divisor.

Throughout the paper, we will assume that $f$ is relatively minimal.
\begin{defn}
We say that a fibred surface is:
\begin{itemize}
\item[-] \textit{smooth}  if every fibre is smooth;
\item[-] \textit{isotrivial} if all smooth fibres are mutually isomorphic;
\item[-] \textit{locally trivial} if $f$ is smooth and isotrivial (equivalently if $f$ is a fibre bundle):
\item[-] {\em trivial} if $S$ is birationally equivalent to $F\times B$ and $f$ corresponds to the projection on $B$. If $b>0$ and $f$ is relatively minimal this is equivalent to $S=F\times B$, 
\end{itemize} 
\end{defn}

Recall the following relative numerical  invariants for fibred surfaces:
\begin{itemize}
\item $K_f^2=K_S^2-8(g-1)(b-1)$ the self-intersection of the relative canonical divisor;
\item $\chi_{f}:=\chi(\mathcal{O}_S)-(g-1)(b-1)= \deg f_{*}\omega_f$ the relative Euler characteristic (the last equality follows from Leray's spectral sequence);
\item $e_{f}:= e(S)-e(B)e(F)=e(S)-4(g-1)(b-1)$ the relative topological characteristic (with $e(X)$ topological characteristic of $X$);
\item $q_{f}:=q-b$ the relative irregularity, with $q=h^1(S,\mathcal{O}_{S})$ irregularity of $S$.
\end{itemize}
For those invariants the following relations are known \cite{Arakelov},\cite{ACGH2}, \cite{Beau}:
\begin{enumerate}
\item $K_{f}^2\geq 0$ and $K_{f}^2=0$ if and only if $f$ is locally trivial (see Remark \ref{rem: lt});
\item \label{ii} $\chi_{f}\geq 0$ and $\chi_f =0$ if and only if $f$ is locally trivial;
\item $e_{f}\geq 0$ and $e_{f}=0$ if and only if $f$ is smooth;
\item $q_f\leq g$ and equality holds if and only if $f$ is trivial.
\end{enumerate}
From Groethendieck-Riemann-Roch theorem we have Noether's relation \cite{ACGH1}
$$ 12\chi_f = K_f^2 +e_f .$$
\begin{rem} \label{rem: lt}
Suppose that $K_f^2=0$. Then by the slope inequality we have $\chi_f=0$ so $f$ is locally trivial. If, on the other hand, $f$ is locally trivial; then  by $(ii)$ we have  $\chi_f=0$, then by Noether's relation and the non-negativity of $e_f$ we have $K_f^2=0$.
\end{rem}
\begin{defn}
The rank $g$ vector bundle $f_{*}\omega_{f}$ is called the \textit{Hodge bundle} of the fibred surface.
\end{defn}
We have the following decompositions of the Hodge bundle as a direct summand of vector sub-bundles:
\begin{itemize}
\item (First Fujita decomposition \cite{Fuj1})
\begin{equation}\label{first fujita}
f_{*}\omega_{f}=\mathcal{O}_{B}^{\oplus q_{f}}\oplus \mathcal{E}, 
\end{equation}
 where $\mathcal{E}$ is nef and  $H^0(B,\mathcal{E}^\vee)=0$;
\item  (Second Fujita decomposition \cite{Fuj2} \cite{CD2})
\begin{equation}\label{second fujita} 
f_{*}\omega_{f}=\mathcal{A}\oplus \mathcal{U},
\end{equation}
with $\mathcal{A}$ ample and $\mathcal{U}$ unitary flat.
\end{itemize}
\begin{defn}\label{unitary rank}
Following \cite{SGT}, we define the \textit{unitary rank} $u_{f}$ of the fibred surface to be the following integer
$$u_{f}:=\rk \mathcal{U}. $$
\end{defn}
\begin{rem}\label{rem: uf}
Comparing the two decomposition, since every trivial bundle is unitary flat, we have:
$$\mathcal{O}_{B}^{\oplus q_{f}}\subseteq \mathcal{U}, $$
and then it holds that $q_{f}\leq u_{f}$. 
Moreover, $\deg \mathcal {U}=0$ and $\deg \mathcal A>0$, hence 
$$\chi_f=\deg f_{*}\omega_f=\deg \mathcal {A},$$ and 
$u_{f}=g$ if and only if $\chi_f=0$ (equivalently $f$ is locally trivial). 
Catanese and Dettweiler first gave examples \cite{CD} \cite{CD2} of fibred surfaces for which the unitary summand is not semiample, thus disproving a long standing conjecture of Fujita. They proved that semi-ampleness of the Hodge bundle is indeed equivalent to $\mathcal{U}$ having finite monodromy. In all the examples in loc. cit. $q_f=0$, hence in particular the strict inequality $q_{f}<u_{f}$ holds. Note moreover that for any fibred surface such that  the monodromy of $\mathcal{U}$ is infinite, the inequality $q_{f}<u_{f}$ also  holds ``up to base change'', i.e.  for any fibration $\tilde f$ obtained from $f$ via base change, we still have $q_{\tilde f}<u_{\tilde f}$. On the other hand, if the monodromy is finite, then there exist a base change $a\colon \tilde B\to B$ such that the induced fibration $\tilde{f}$ has $q_{\tilde{f}}=u_{\tilde{f}}$. See \cite{SGT}. \end{rem}

Over the moduli space $\mathcal{M}_{g}$ of smooth curves of genus $g$, the function:
$$ [C] \mapsto \cliff(C) $$
is a well defined lower semicontinuous function 
This allows us to give the following:

\begin{defn}
Given $f\colon S\rightarrow B$ a relatively minimal fibred surface. We define 
$$c_{f}:= \max_{t \in B}\{\cliff(F_{t}) \ | \text{$F_{t}$ is a smooth fibre of $f$} \}=\cliff(F)\mbox{ for $F$ general fibre of $f$}$$
and call it the \textit{Clifford index of f}.
\end{defn}

\subsection{Xiao's technique}\label{subsec: xiao}

In this section we recall the main results of Xiao's method, introduced by Xiao in his seminal paper \cite{Xiao}, and further developed by Konno and Barja. We will then apply this method to a subbundle of the Hodge bundle, but we think it is worth to develop in full generality and detail the construction, as the precise statement we need is not immediate to find in the literature.
Let  $ \pi \colon \mathbb{P}_{B}(\mathcal{E})\rightarrow B$
be the projective bundle of one dimensional quotients of $\mathcal{E}$ (Grothendieck's notations);
and let $\mathcal{O}_{\mathbb{P}(\mathcal{E})}(1)$ be the  associated tautological line bundle.
\begin{defn}
We say that $\mathcal{E}$ is a \textit{nef (resp. ample) vector bundle} if   $\mathcal{O}_{\mathbb{P}(\mathcal{E})}(1)$ is nef (resp. ample) over $\mathbb{P}_{B}(\mathcal{E})$.
\end{defn}
Let $f\colon S\rightarrow B$ be a relatively minimal fibration and fix a divisor $D$ on $S$.
For every non zero vector subbundle $\mathcal{F}\subseteq f_{*}\mathcal{O}_S(D)$, the natural homomorphism $$f^{*}\mathcal{F}\hookrightarrow f^{*}f_{*}\mathcal{O}_S(D)\longrightarrow \mathcal O_S(D)$$ yields a rational map 
$$
\xymatrix{
S\ar@{-->}^{\psi}[rr]\ar_f[dr] &&\mathbb{P}_{B}(\mathcal{F})\ar^{\pi}[dl]\\
&B&}
$$
such that $\pi \circ\psi = f$.
The indeterminacy locus of the map $\psi$ is described by the following result, whose proof is immediate.
\begin{thm}\label{ohno}[Ohno \cite{Onno}] In the above situation, there exists a blow up $\epsilon\colon \hat{S} \rightarrow S$ and a morphism $\lambda :=\psi \circ \epsilon \colon \hat{S}\rightarrow \mathbb{P}_{B}(\mathcal{F})$
such that
$\lambda^{*}L_{\mathcal{F}} \sim \epsilon^{*}(D-Z)-E$ where 
\begin{itemize}
\item $Z$ is an effective divisor on $S$;
\item $E$ is a $\epsilon-$exceptional effective divisor of $\hat{S}$;
\item  $L_{\mathcal{F}}$ a  hyperplane section of $\mathbb{P}_{B}(\mathcal{F})$ i.e. a divisor associated to $\mathcal{O}_{\mathbb{P}(\mathcal{F})}(1)$. 
\end{itemize}
\end{thm}     
\begin{defn}
In this setting we define:
\begin{itemize}
\item $M(D,\mathcal{F}):= \lambda^{*}L_{\mathcal{F}}$  the \textit{moving part} of the vector subbundle $\mathcal{F}$;
\item $Z(D,\mathcal{F}):=\epsilon^{*}Z +E$  the \textit{fixed part} of the vector subbundle $\mathcal{F}$;
\item $N(D,\mathcal{F}):= M(D,\mathcal{F})-\lambda^{*}\mu(\mathcal{F})F$ where we note that $\epsilon$ do not change the general fibre of $f$; then we can rewrite:
$N(D,\mathcal{F})=M(D,\mathcal{F})-\mu(\mathcal{F})F$ with $F$ a general fibre of $f$.
\end{itemize}
\end{defn}
The Xiao's method makes a crucial use of the Harder-Narasimhan filtration.
\begin{defn}
Let $\mathcal{F}$ a vector bundle over a smooth projective curve $B$. There exists a unique sequence of vector sub-bundles of $\mathcal{F}$:
\begin{equation}
0= \mathcal{F}_{0}\varsubsetneq \mathcal{F}_{1}\varsubsetneq \ldots\varsubsetneq \mathcal{F}_{k-1}\varsubsetneq \mathcal{F}_{k}=\mathcal{F}
\label{A}
\end{equation}
satisfying the conditions:
\begin{itemize}
\item for $i=1,\ldots,k$ $\mathcal{F}_{i}/\mathcal{F}_{i-1}$ is a semistable vector bundle;
\item For any $i=1,\ldots,k$ setting $\mu_{i}:=\mu(\mathcal{F}_{i}/\mathcal{F}_{i-1})$, we have:
$$\mu_{1}> \mu_{2} >\ldots>\mu_{k}. $$
\end{itemize}
The filtration \ref{A} is called \textit{Harder-Narasimhan filtration} of $\mathcal{F}$.

We set $\mu_{-}(\mathcal{F}):=\mu_{k}$,  and call  it the {\em final slope} of the sheaf. 
\end{defn}
\begin{rem}
Note that it holds the formula:
$$\deg\mathcal{F}= \sum_{i=1}^{k}{r_{i}(\mu_{i}-\mu_{i+1})}.$$
Indeed, 
considering the exact sequence of vector bundles:
$$0\rightarrow \mathcal{F}_{k-1}\rightarrow \mathcal{F}_{k}\rightarrow \mathcal{F}_{k}/\mathcal{F}_{k-1}\rightarrow 0, $$
from the additivity property of degree, we can say
$\deg \mathcal{F}_{k}= \deg \mathcal{F}_{k-1}+ \deg \mathcal{F}_{k}/\mathcal{F}_{k-1}. $
Similarly,  we have that:
$\deg \mathcal{F}_{k-1}= \deg \mathcal{F}_{k-2}+ \deg \mathcal{F}_{k-1}/\mathcal{F}_{k-2},$
and so on. By induction we can conclude that:
$$\deg \mathcal{F}_{k}= \deg(\mathcal{F}_{k}/\mathcal{F}_{k-1}) +\deg(\mathcal{F}_{k-1}/\mathcal{F}_{k-2})+ \ldots+\deg(\mathcal{F}_{2}/\mathcal{F}_{1})+\deg (\mathcal{F}_{1}) = \sum_{i=1}^{k}{\deg (\mathcal{F}_{i}/\mathcal{F}_{i-1})}.$$
Now, from the definition of slope, for every $i=1,..,k$ we have $\deg \mathcal{F}_{i}/\mathcal{F}_{i-1}= \mu_{i}(r_{i}-r_{i-1}),$
So, setting $\mu_{k+1}=0$ and $r_{k+1}=r_k$, we obtain the desired formula
$$\deg\mathcal{F}=\deg \mathcal{F}_{k}= \sum_{i=1}^{k}{\mu_{i}(r_{i}-r_{i-1})}= \sum_{i=1}^{k}{r_{i}(\mu_{i}-\mu_{i+1})}.$$
\end{rem}
The Xiao's method is based on the following fundamental result of Miyaoka-Nakayama.
\begin{thm}(\cite{Miy} 
\cite[Corollary 3.8]{nak})
[Miyaoka-Nakayama]
\label{Miy} Let $\mathcal{F}$ be a locally free sheaf on a projective curve $B$.
Let $\Sigma$ be the general fibre of $\pi :\mathbb{P}_{C}(\mathcal{F})\rightarrow C$. The $\mathbb{Q}$-divisor $L_{\mathcal{F}}-x\Sigma$ is nef if and only if $x\leq \mu_{-}(\mathcal{F})$. 
\end{thm}
\begin{rem}\label{rem: nef}
From Miyaoka-Nakayama's result we see straightforwardly that $\mu_{-}(\mathcal{F})\geq 0$ if and only if $\mathcal{F}$ is a nef vector bundle on $B$.
\end{rem}
\begin{rem}\label{l-1}
In the case $\mathcal G=f_*\omega_f$, it is important to notice that the second to last subsheaf is precisely the ample part in second Fujita's decomposition:
$\mathcal{F}_{l-1}= \mathcal{A}.$
Indeed, $f_*\omega_f$ is nef, and the subsheaf $\mathcal{U}=f_*\omega_f/\mathcal{A}$ is a subsheaf of maximal rank in $f_*\omega_f$ with (minimal) degree 0.
For the Hodge bundle the last slope $\mu_l$ is greater or equal to $0$ and  $\mu_l=0$ if and only if $\mathcal U\not =0$.
\end{rem}
We are now ready to expose the heart of Xiao's method:
\begin{thm}(Xiao's key Lemma \cite{Xiao})
Let $f\colon S\rightarrow B$ be a fibred surface. Let $D$ be a divisor on $S$ and suppose that there exist a sequence of effective divisors:
$$Z_{1}\geq Z_{2}\geq \ldots \geq Z_{s}\geq Z_{s+1}:=0 $$
and a sequence of rational numbers 
$$\mu_{1} > \mu_{2} > \ldots\ldots > \mu_{s}\geq \mu_{s+1}:=0 $$
such that for every $i=1,\ldots,s$ $N_{i}:=D-Z_{i}-\mu_{i}F$ is a nef $\mathbb{Q}-$divisor. Then for any set of indexes $\{j_1,\ldots ,j_s\}\subseteq \{1,\ldots ,\l\}$ we have
$$ D^{2}\geq \sum_{i=1}^{s}{(d_{j_i}+d_{j_{i+1}})(\mu_{j_i}-\mu_{j_{i+1}})} $$
where $d_{j}:=N_{j}F$.
\label{xiao}
\end{thm}

\begin{proof}
Just observe that the assumptions imply the following:
\begin{equation*}
\begin{split}
N_{j_{i+1}}^2-N_{j_i}^2=&\,\,\,(N_{j_{i+1}}+N_{j_i})(N_{j_{i+1}}-N_{j_i})= (N_{j_{i+1}}+N_{j_i})(Z_{j_{i}}-Z_{j_{i+1}}-(\mu_{i}-\mu_{i+1})F)\\
\geq &\,\,\,(d_{j_{i}}+d_{j_{i+1}})(\mu_{i}-\mu_{i+1}),\\
\end{split}
\end{equation*}
and that 
$$\sum_{i=1}^s(N_{j_{i+1}}^2-N_{j_i}^2)=-N_{j_{1}}^2+N_{j_s}^2\leq N_{j_s}^2\leq D^2.$$
\end{proof}

\subsection{Main inequality}\label{subsec: mainxiao}
We are now ready to state the version of Xiao's basic result in the form needed. Note that this is an expanded version of the inequality stated in \cite[Remark 24]{SB-survey}.
\begin{thm} \label{relative}
Let $f\colon S\rightarrow B$ be a fibred surface. Let $D$ be a nef divisor on $S$ and $\mathcal{G}\subseteq f_{*}\mathcal{O}_{S}(D)$ be a rank $r$ subsheaf. 
Let $d'= MF$ where $M=M(D,\mathcal G)$.

Suppose that there exists a real number  $\alpha>0$ such that  for every linear subsystem $|P|$ of $|M_{|F}|$
\begin{equation}\label{ass-alpha}
\frac{\deg |P|}{\dim |P|}\geq \alpha. 
\end{equation}
\begin{enumerate}
\item The following inequality holds:
\begin{equation}\label{rel2}
D^2\geq \frac{2\alpha(r-1)}{r}\deg\mathcal{G}=2\alpha (r-1)\mu( \mathcal{G}). 
\end{equation}
\item If moreover $\mathcal G$ is nef, then, for every  non negative integer $d\leq d'$,  the following inequality holds:
\begin{equation}
D^2\geq \frac{2\alpha d}{d+\alpha}\deg \mathcal{G}.
\end{equation}
\end{enumerate}
\end{thm}
\begin{proof}
Let
\begin{equation}
 0 \varsubsetneq \mathcal{G}_{1}\varsubsetneq\ldots.\varsubsetneq \mathcal{G}_{k-1}\varsubsetneq \mathcal{G}_{k}=\mathcal{G}
\label{scomp} 
\end{equation}
be the Harder-Narasimhan filtration of $\mathcal{G}$. We note that in general this filtration need not necessarily be related to the Harder-Narasimhan filtration of $f_{*}\mathcal{O}_{S}(D)$ (although this will happen in the application: see Remark \ref{restrizione HN}). 

Following Ohno's construction in Theorem \ref{ohno}, we consider  a suitable blow up $\nu\colon \hat{S}\to S$ and over $\hat S$ for every index $i$ we consider  the divisors $M_i:=M(D,\mathcal{G}_{i})$ and $Z_i:=Z(D,\mathcal{G}_{i})$, which are respectively nef and effective.
Call $r_{i}=\rk\mathcal{G}_{i} $ and $d_{i}:=M_iF$. 
 We also set 
$\mathcal{G}_{k+1}:=\mathcal{G}_{k}=\mathcal{G}.$ 

%

Let us first assume that $\mathcal G$ is nef and prove inequality (\ref{relative}). The final slope of $\mathcal G$ is  $\mu_{k}\geq 0$ by Remark \ref{rem: nef} and we can choose $\mu_{k+1}=0$ and $ Z_{k} = Z_{k+1}$.
The sequence $(Z_{i},\mu_{i})$ clearly satisfies by construction: 
$$ Z_{1}\geq Z_{2}\geq\ldots\geq Z_{k} = Z_{k+1}, $$ 
and 
$$ \mu_{1}> \mu_{2}>\ldots>\mu_{k}\geq \mu_{k+1}:=0.$$ 
Observing that $\mu_{i} $ coincides with $\mu_{-}(\mathcal{G}_{i})$ 
we have by the Theorem of Miyaoka-Nakayama  that the divisors
$$N_{i}:= M(D, \mathcal{G}_{i})-\mu_{i}F$$
are all nef $\mathbb{Q}-$divisors over $\hat{S}$.
Since the intersection product is invariant under birational morphism we have $(\nu^{*}D)^2=D^2. $
So, we can apply Theorem  \ref{xiao} to estimate 
 $(\nu^{*}D-Z_{k})^2$. We make a wise use of the choice of the indexes in the theorem.
 
Firstly we use the set of indexes $\{ 1,\ldots,k\}$, obtaining the inequality
$$(\nu^{*}D-Z_{k})^{2}\geq \sum_{i=1}^{k}{(d_{i}+d_{i+1})(\mu_{i}-\mu_{i+1})}, $$
which in its extensive form reads as follows
$$(\nu^{*}D-Z_{k})^2\geq (d_{1}+d_{2})(\mu_{1}-\mu_{2}) +\ldots+(d_{k-1}+d_{k})(\mu_{k-1}-\mu_{k})+(d_{k}+d_{k+1})(\mu_{k}).$$
Observe that  assumption (\ref{ass-alpha}) implies that for any $i$, $d_i\geq \alpha(r_i-1)$, because in case $r_1=1$, the inequality holds trivially.
Using this inequality and the fact that $r_{i}\geq r_{i-1}+1$ for $i=1,\ldots,k-1$ and that $r_{k+1}=r_{k}$, we have:
\begin{equation*}
\begin{split}
(\nu^{*}D-Z_{k})^{2}  &   \geq  \sum_{i=1}^{k}{(d_{i}+d_{i+1})(\mu_{i}-\mu_{i+1})} \geq \\       
       &    \geq 2\alpha(\sum_{i=1}^{k-1}{r_{i}(\mu_{i}-\mu_{i+1})}+r_{k}\mu_{k}) -\alpha(\mu_{1}+\mu_{k}) =  \\
          & = 2\alpha \deg \mathcal{G} -\alpha(\mu_{1}+\mu_{k}).
\end{split}
\end{equation*}
Consider now the list of indexes $\{ 1,k\}$: we have
$$(\nu^{*}D-Z_{k})^{2}\geq (d_{1}+d_{k})(\mu_{1}-\mu_{k})+(d_{k}+d_{k+1})(\mu_{k})\geq d_{k}(\mu_{1}+\mu_{k}). $$ 
Eventually, combining the last two inequalities we obtain:
$$ (\nu^{*}D-Z_{k})^{2}\geq \frac{2\alpha d_{k}}{d_{k}+\alpha}\deg \mathcal{G}.$$
Now observe that 
$$(\nu^{*}D-Z_{k})^2= D^2 -2\nu^{*}DZ_{k}+Z_{k}^2\leq D^2,$$ 
where the last inequality follows from the fact that $\nu^{*}D$ is nef and  $Z_{k}$ effective and from $Z_{k}^2\leq 0$ by Hodge index theorem.
Now, consider the following function:
$$h(t):=\frac{2\alpha t}{\alpha + t}, $$
which is monotonically increasing for $t\geq 0$.
From the hypothesis we have
$d_{k}\geq d $
so we can deduce that 
$$D^2\geq \frac{2\alpha d_{k}}{d_{k}+\alpha}\deg \mathcal{G} = h(d_k)\deg \mathcal{G}\geq h(d)\deg \mathcal{G}=\frac{2\alpha d}{d+\alpha}\deg \mathcal{G},$$
and the proof of inequality (\ref{relative}) is concluded under the assumption that $\mathcal G$ is nef. 

In the non-nef case, just consider as in \cite[Prop.8]{SB-survey} the last nef subbundle in the Harder-Narasimhan sequence: 
$\mathcal G_s$, where $s=\max\{i\mid \mu_i\geq 0\}$. Applying the very same construction to $\mathcal G_s$ we can obtain 
$$D^2\geq\frac{2\alpha d_s}{d_s+\alpha}\deg \mathcal{G}_s\geq 2\frac{\alpha(r_s-1)}{r_s}\deg\mathcal{G}_s\geq 2\frac{\alpha(r-1)}{r}\deg\mathcal{G},$$
where the second inequality is obtained by choosing $d=\alpha(r_s-1)$, and the last inequality follows from the monotonicity of the function  $h(t)$ above and from the fact that clearly $\deg{\mathcal G}_s\geq \deg \mathcal G$.
So, also inequality (\ref{rel2}) is proved.
\end{proof}
\begin{rem}
As proved in \cite{barja-ine}, the vector subbundle $\mathcal G_s$ in the proof of the above theorem is a maximal element in the set of nef sub-bundles of $\mathcal G$: for any nef sub-bundle of $\mathcal G$ it holds $\mathcal F\subseteq \mathcal G_s$.
\end{rem}

In particular, if $|\mathcal{G}\otimes \C(t)|$ is linearly semistable for general $t\in B$,  we can take:
$$\alpha = \frac{\deg |\mathcal{G}\otimes \C(t)|}{\dim |\mathcal{G}\otimes \C(t)|}.$$
and obtain the following well known result (see \cite{SB-survey}).
\begin{cor}
Let $f\colon S\rightarrow B$ be a fibred surface. Given $D$ a nef divisor on $S$ and $\mathcal{G}\subseteq f_{*}\mathcal{O}_{S}(D)$ a rank $r$ subsheaf. Let  $d=\deg |\mathcal{G}\otimes \C(t)|$ the degree of the linear system $|\mathcal{G}\otimes \C(t)|$, over a general fibre $F_{t}$.
If $|\mathcal{G}\otimes \C(t)|$ is linearly semistable, then
$$D^2\geq \frac{2d}{r}\deg \mathcal{G}= 2d\mu(\mathcal{G}).$$
\end{cor}

\section{Slope inequalities}\label{sec: main}

Let $f\colon S\rightarrow B$ be a relative minimal fibration of genus $g\geq 2$. 
We are now ready to prove our main estimates on the slope of fibred surfaces. 

Firstly, using the first Fujita decomposition (\ref{first fujita}) we give a bound that improves the main bound of \cite{SB}.
Note that the proof is much simpler than the proof of \cite{SB}, where we needed to lift a general projection on the fibre to obtain the desired subsheaf of the Hodge bundle.
 \begin{thm}\label{sim} Let  $m:=\min\{ q_{f},c_{f}\}$. The following inequality holds:
$$K_{f}^{2}\geq 2\frac{2g-2-m}{g-m}\chi_{f}. $$ 
\end{thm}
\begin{proof}
First observe that in the hyperelliptic case $m=0$ and  the inequalities are just the classical slope inequality. Assume that the general fibre is not hyperelliptic.

Let us consider the  first Fujita decomposition (\ref{first fujita}).
$$f_{*}\omega_{f}= \mathcal{E}\oplus \mathcal{O}^{\oplus q_f}.$$
If $q_{f}\leq c_{f}$ consider  the vector bundle $\mathcal{G}:=\mathcal{E}.$
If $q_{f}\geq c_{f}$ consider the vector bundle $\mathcal{G}:=\mathcal{E}\oplus \mathcal{O}_{B}^{q_{f}-c_{f}}$
In both cases the fibre over a general $t\in B$
$\mathcal{G}\otimes \C(t)\subseteq H^0(F_{t},K_{F_{t}}) $
defines a linear subsystem of $H^0(F_{t},K_{F_{t}})$ of codimension $m$.

Let us start by observing that in case that the first vector subbundle in the Harder-Narasimhan filtration of the Hodge bundle is of rank one (a line bundle), we have $d_1=0=r_1-1$. 
By the remark above and Theorem \ref{Lemma1}, we can apply Theorem \ref{relative} to $D=K_{f}$ and $\mathcal{G}$ as defined above, with $\alpha=\frac{2g-2-m}{g-m-1}$. 
We thus obtain
$$K_{f}^2\geq \frac{2\alpha d}{\alpha +d}\deg\mathcal{G}= 2\frac{2g-2-m}{g-m}\chi_f ,$$
as desired.
\end{proof}
We shall now turn our attention on the influence of the unitary rank $u_f$ on the slope. 

\begin{thm}\label{Y}
The following inequalities holds:
\begin{equation*}
K_{f}^{2} \geq \left\{ \begin{array}{ll}
                                      2\frac{2g-2-u_{f}}{(g-u_{f})}\chi_{f} \ \ \ \ \ \ \ \ \ \ \ \ \text{$u_{f}\leq c_{f}$}  \\
                                      2\frac{(2g-2-c_{f})(g-u_{f}-1)}{(g-c_{f}-1)(g-u_{f})}\chi_{f} \ \ \ \text{$u_{f}\geq c_{f}$}
                                     \end{array} \right.  
\end{equation*}

\end{thm}

\begin{proof}
As above, we can assume thet $F$ is non-hyperelliptic.
Consider the second Fujita decomposition (\ref{second fujita})
$f_{*}\omega_{f}= \mathcal{A}\oplus \mathcal{U}$.
As already observed, we have that $\deg \mathcal{A}=\deg f_{*}\omega_{f}$.
We distinguish the two following cases:
\begin{itemize}
\item If $u_{f}\leq c_{f}$, then consider  $\mathcal G=\mathcal{A}$.
From Theorem \ref{A} we can estimate the degree of that linear subsystem as follows:
$$\deg |\mathcal{A}\otimes \C(t)|\geq \frac{2g-2-m}{g-m-1}(g-u_{f}-1)=2g-2-u_{f}=:d.$$
Then, applying Theorem \ref{relative} with $D=K_{f}$ and $\mathcal{G}=\mathcal{A}$, we have:
$$K_{f}^2\geq \frac{2\alpha d}{\alpha +d}\deg \mathcal{A}= 2\frac{2g-2-u_{f}}{g-u_{f}}\chi_{f},$$
as wanted.

\item If $u_{f}\geq c_{f}$, using Theorem \ref{Lemma1} we estimate the degree of the linear system $|\mathcal{A}\otimes\C(t)|$ as:
$$\deg |\mathcal{A}\otimes \C(t)|\geq \frac{2g-2-c_{f}}{g-c_{f}-1}(g-u_{f}-1)=:d. $$ 
Then applying Theorem \ref{relative} with $D=K_{f}$ and $\mathcal{G}=\mathcal{A}$ we have:
$$K_{f}^2\geq \frac{2\alpha d}{\alpha +d}\deg \mathcal{A}= 2\frac{(2g-2-c_{f})(g-u_{f}-1)}{(g-u_{f})(g-c_{f}-1)}\chi_{f},$$
and the proof is concluded.
\end{itemize}
\end{proof}


\begin{rem}
Observe that these last  inequalities are not symmetric in $\min\{u_f,c_f\}$ as  the one of Theorem \ref{sim}. 
In case there exists a unitary flat subsheaf $\mathcal{U'}$ of $\mathcal{U}$, with $\rk \mathcal{U'}\geq  u_{f}-c_{f}$, one can improve the last inequality of Theorem \ref{Y}. 
However, such a subsheaf $\mathcal U'$ need not to exist.  
\end{rem}

\begin{rem}\label{restrizione HN}
It is worth making the following remark. In Xiao's method as exposed in Section \ref{subsec: xiao}, we use the Harder-Narasimhan sequence of the subsheaf $\mathcal{G}$ of $f_*\mathcal{O}_S(D)$. This in general is not related to the Harder-Narasimhan sequence of $f_*\mathcal{O}_S(D)$ itself. 
But in case $\mathcal{G}$ is a nef subsheaf of the Hodge bundle that contains the ample summand $\mathcal{A}$, then the Harder-Narasimhan filtration of $\mathcal{G}$ clearly is the truncation of the filtration of $f_{*}\omega_{f}$. 

\end{rem}
\section{An example}\label{example}

Now we want to expose a first example of a fibred surface in which the bound of Theorem \ref{Y} is better than the bound of Theorem \ref{sim} and of any other previous bound. The known examples of fibred surfaces with high unitary rank (\cite{CD,CD2,CD3, Lu}) all satisfy  $c_f\leq 1$. It is therefore  interesting to find examples with Clifford index close to the unitary rank. This is a first example in this direction. 

We use the same construction of \cite{CD3}, although with some modifications, and refer to loc. cit. for details.
Consider a family of cyclic coverings of $\mathbb{P}^1$ of degree 7, given birationally by the equation:
\begin{equation}\label{curva}
y^7= x_{0}x_{1}(x_{1}-x_{0})^{2}(x_{1}-tx_{0})^{3},
\end{equation}
where $x_{0},x_{1}$ are homogeneous coordinates of $\mathbb{P}^1$, $t\in \C \setminus \{0,1 \}$ and we can think $y$ as a section of $\mathcal{O}_{\mathbb{P}^1 \times \mathbb{P}^1}(1,1)$.\
   
By base change and compactification we can associate a non-locally trivial fibred  surface $S\colon \to B$ with the following properties:
\begin{itemize}
\item the base curve $B$ is of genus $b=3$;
\item the general fibre $F$ is a smooth curve of genus $g=6$;
\item there are only three singular fibres, which are given by two smooth curves of genus 3 intersecting transversely (in particular $f$ is a semistable fibration); 
\item on every fibre there is an action of $\Z/7\Z$ by automorphisms;
\item the irregularity $q=h^1(S,\mathcal{O}_{S})$ is equal to 3.
\end{itemize}  

This last property allows us to say that $f\colon S\rightarrow B$ is an Albanese fibration, i.e. that $q_f=0$.

The  action  of $\Z/7\Z$ on the smooth fibres  the space induces an action on $H^{1}(F,\mathcal{O}_{F})=H^{0}(F,\omega_{F})$ and for any character $\chi_{j}\in \Z/7\Z$ we can calculate  the dimension of the corresponding characteristic subspace $H^{0}(F,\omega_{F})_{j}$ via the Chevalley-Weil formula:
\begin{equation*}
\dim H^{0}(F, \omega_{F})_{j}= -1 +\sum_{k=1}^{4}{\frac{(-j\alpha_{k})_{7}}{7}} = \left\lbrace  \begin{array}{l}
       2 \ \ \ \ \ \text{if $j=1$} \\
       1 \ \ \ \ \ \text{if $j=2,3,4,5$} \\
       0 \ \ \ \ \ \text{if $j=6$.}
\end{array} \right.
\end{equation*}

where $(\alpha_{1}, \alpha_{2}, \alpha_{3},\alpha_{4})=(1,1,2,3)$ and $(n)_{7}$ means $n \mod 7$.

Associated to these subspaces we have sub-bundles of the Hodge bundle, and by the same arguments of \cite{CD3} one can see that the rank 2 sub-bundle associated to $j=1$ is contained in $\mathcal U$. 
On the other hand, for $j=2,3,4,5$ we have rank one summands, hence line bundles. If these line bundles were contained in the unitary part, they necessarily would have infinite monodromy, which is impossible, as proved by Deligne (see \cite[Corollary 21]{CD3}). So, we have $u_f=2$ for this fibration.

Moreover, since the general fibre $F$ has equation \ref{curva}, applying the results in \cite{gon1} and \cite{gon2}, we conclude that $F$ has maximal gonality 4.
The Clifford index of $F$ therefore is either 2 or 1 in case $F$ possesses a $g^5_2$. We want to exclude this last case: observe that $F$ has an automorphism of order 7,  while for a plane curve of degree 5 the order of any cyclic subgroup of the automorphism has to divide one of the following integeres
\[4,\,5,\,10,\,16,\,15,\,20\]
by \cite[Cor.8]{BB}.
Observe moreover that the general fibre is not bielliptic as can be derived for instance by \cite[Lemma 2.4]{schweizer}, so the results of \cite{barja-biell} do not apply.

Now, from  Theorem \ref{Y}, since $u_{f}=c_{f}=2$ and $g=6$, 
\begin{equation*}
K_{f}^2/\chi_{f}\geq \frac{2(2g-2-u_{f})}{g-u_{f}}= 4.
\end{equation*}  
By direct computations, one can see in this case that we have $K_{f}^2/\chi_{f}= 45/4,$
that is strictly greater than 4.

It remains open the question wether the inequalities proved in this paper are sharp.

\bigskip

\noindent Dottorato Milano Bicocca-Pavia CICLO XXXIV,
 \\Dipartimento di Matematica,Universit\`a di Pavia,\\ \textsl {e.riva55@campus.unimib.it}.
 
\medskip

\noindent Dipartimento di Matematica, Universit\`a di Pavia,\\
 \textsl {lidia.stoppino@unipv.it}.


\begin{thebibliography}{ABCD}


\bibitem{ACGH1}  E.~Arbarello, M.~Cornalba, P.A.~Griffiths, J.~Harris, {Geometry of algebraic curves. {V}ol. {I}}, {GMW}, {\bf 267}, {Springer-Verlag, New York}, {1985}.

\bibitem{ACGH2}  E.~Arbarello, M.~Cornalba, P.A.~Griffiths, {Geometry of algebraic curves. {V}ol. {2}}, {GMW}, {\bf 267}, {Springer-Verlag, New York}, {2010}.

\bibitem{Arakelov} S.J.~Arakelov, {\em Families of algebraic curves with fixed degeneracies,} Math. USSR Izvestija, {\bf 5 } no. 6 (1971),1269--1293.

\bibitem{BB} E.~Badr, F.~Bars {\em Plane non-singular curves with an element of ``large? order in its automorphism group}, Int. Journal of Algebra and Computation {\bf 26}, no. 2 (2016), 399--433.

\bibitem{Ballico} E.~Ballico, {\em On the Clifford index of algebraic curves}, Proc. Amer. Math. Soc. {\bf 97} (1986), 217--218.


\bibitem{barja-biell} M.A.~Barja, {\em On the slope of bielliptic fibrations}. Proc. Amer. Math. Soc. {\bf129} (2001), no. 7, 1899--1906. 

\bibitem{barja-ine} M.A.~Barja, {\em Higher dimensional slope inequalities for irregular fibrations}, arXiv:2012.06889 [math.AG].


\bibitem{GNB} M.A.~Barja, V.~Gonz\`alez-Alonso, J.C.~Naranjo, {\em Xiao's conjecture for general fibred surfaces}, J. reine angew. Math. {\bf 739} (2018), 297--308.

\bibitem{SB} M.A.~Barja, L.~Stoppino, {\em Linear stability of projected canonical curves with applications to the slope of fibred surfaces},  J. Math. Soc. Japan, {\bf 60}, No. 1 (2008), 171--192.

\bibitem{BS-trigonal} M.A.~Barja, L.~Stoppino, {\em Slopes of trigonal fibred surfaces and of higher dimensional   fibrations}, {Ann. Sc. Norm. Super. Pisa Cl. Sci. (5)}, {\bf 8}, (2009), {No. 4}, { 647--658}.

\bibitem{SB-survey} M.A.~Barja, L.~Stoppino, {\em Stability conditions and positivity of invariants of fibrations}, {Algebraic and complex geometry}, {Springer Proc. Math. Stat.}, {\bf 71}, {1--40}, {(2014)}.


\bibitem{Beau} A.~Beauville, Appendix to {\em In\'egalit\'es num\'eriques pour les surfaces de type g\'en\'eral.}, by O. Debarre, Bull. Soc. Math. France {\bf 110} (1982), no. 3, 319--346.
 
\bibitem{BZ} V.~Beorchia, F.~Zucconi, {\em On the slope of fourgonal semistable fibrations},  Math. Res. Lett. {\bf 25}, no. 3 (2018), 723--757. 

\bibitem{CD} F.~Catanese, M.~Dettweiler, {\em The direct image of the relative dualizing sheaf needs not be semiample}, C. R. Math. Acad. Sci. Paris {\bf 352} (2014), no. 3, 241--244.

\bibitem{CD2} F.~Catanese, M.~Dettweiler, {\em Vector bundles on curves coming from variation of Hodge structures}, Internat. J. Math. {\bf 27} (2016), no. 7, 1640001, 25 pp.

\bibitem{CD3} F.~Catanese, M.~Dettweiler, {\em Answer to a question by Fujita on variation of Hodge structures}, Higher Dimensional Algebraic Geometry (Tokyo: Mathematical Society of Japan, 2017), 73--102.

\bibitem{CLZ}   K.~Chen, X.~Lu, K.~Zuo, {\em On the Oort conjecture for Shimura varieties of unitary and orthogonal types}. Compos. Math. {\bf 152 }(2016), no. 5, 889--917. 


\bibitem{Martens} M.~Coppens, G.~Martens,  {\em Secant spaces and Clifford's theorem}, Compositio Math. {\bf 78} (1991), 193--212.

\bibitem{C-H} M.~Cornalba, J.~Harris, {\em Divisor classes associated to families of special varieties, with applications to the moduli space of curves} Ann. Sc. Ec. Norm. Sup., {\bf 21} (1988), no. 4, 455--475.

\bibitem{C-S} M.~Cornalba, L.~Stoppino, {\em A sharp bound for the slope of double cover fibrations}, Michigan Math. J. {\bf 56} (2008), no. 3, 551--561. 

\bibitem{Eno} M.~Enokizono, {\em Slopes of fibered surfaces with a finite cyclic automorphism}, Michigan Math. J. {\bf 66} (2017), no. 1, 125--154. 
 
\bibitem{Fuj1} T.~Fujita, {\em On K\"ahler fiber spaces over curves}, J. Math. Soc. Japan {\bf 30} (1978), no. 4, 779--794.

\bibitem{Fuj2} T.~Fujita, {\em The sheaf of relative canonical form of a K\"ahler fiber space over a curve},  Proc. Japan Acad. Ser. A Math. Sci. {\bf 54} (1978), no. 7, 183--184. 

\bibitem{GKM} A.~Gibney, S.~Keel, I.~Morrison {\em Towards the ample cone of $\overline{\mathcal{M}}_{g,n}$}, J. Amer. Math. Soc. {\bf 15} (2002), no. 2, 273--294.

\bibitem{TG} V.~Gonz\`alez-Alonso, S.~Torelli, {\em Families of curves with Higgs field of arbitrarily large kernel}, to appear in Bullettin of the L.M.S., 2020. DOI: 10.1112/blms.12437

\bibitem{SGT} V.~Gonz\`alez-Alonso, L.~Stoppino, S.~Torelli, {\em On the rank of the flat unitary summand of the Hodge bundle}, Trans. Amer. Math. Soc.{\bf 372} (2019), 8663--8677.




\bibitem{konno-nonhyp} K.~Konno, {\em Nonhyperelliptic fibrations of small genus and certain irregular canonical surfaces.} Ann. Scuola Norm. Sup. Pisa Cl. Sci. (4) 20 (1993), no. 4, 575--595. 

\bibitem{konno-trig} K.~Konno, {\em A lower bound of the slope of trigonal fibrations,} Internat. J. Math. {\bf 7 } (1996), no.1, 19--27.

\bibitem{konno-cliff}  K.~Konno, {\em Clifford index and the slope of fibered surfaces,} J. Algebraic Geom. {\bf 8 } (1999), no. 2, 207--220. 

\bibitem{LZ2} X.~Lu, K.~Zuo, {\em On the slope conjecture of Barja and Stoppino for fibred surfaces,} Ann. Sc. Norm. Super. Pisa Cl. Sci. (5) {\bf 19} (2019), no. 3, 1025--1064.

\bibitem{LZ3} X.~Lu, K.~Zuo, {\em On the gonality and the slope of a fibered surface,}  Adv. Math. {\bf 324} (2018), 336--354. 

\bibitem{Lu} X.~Lu, {\em Family of curves with large unitary summand in the Hodge bundle,} Math. Z. {\bf 291} (2019), no. 3-4, 1381--1387.


\bibitem{MS} {E.C.~Mistretta, L.~Stoppino}, {\em Linear series on curves: stability and {C}lifford index}, {Internat. J. Math.}, {\bf 23 (12)} (2012),{1250121, 25 pp.}

\bibitem{Miy} Y.~Miyaoka, {\em The Chern class and Kodaira dimension of a minimal variety}, Algebraic Geometry, Sendai (1985), 449--476.

\bibitem{moriwaki} A.~Moriwaki,  {\em A sharp slope inequality for general stable fibrations of curves,} J. Reine Angew. Math. {\bf 480} (1996), 177--195.

\bibitem{Mumford} D.~Mumford,  {\em Stability of projective varieties}, Enseign. Math. (2) {\bf 23} (1977), no. 1-2, 39--110. 

\bibitem{nak} N.~Nakayama,  {\em Zariski-decomposition and abundance.} MSJ Memoirs, 14. Mathematical Society of Japan, Tokyo, 2004. xiv+277 pp. 

\bibitem{Onno} K.~Ohno, {\em Some inequality for minimal fibration of surface of general type over curves}, J. Math. Soc. Japan {\bf 44} (1992), no.4, 643--666.

\bibitem{pardini} R.~Pardini, {\em The Severi inequality $K^2\geq4\chi$ for surfaces of maximal Albanese dimension}. Invent. Math. {\bf 159} (2005), no.3, 669--672. 

\bibitem{pirola} G.P.~Pirola, {\em On a conjecture of Xiao.} J. Reine Angew. Math. {\bf 431} (1992), 75--89.


\bibitem{stoppino} L.~Stoppino, {\em Slope inequalities for fibred surfaces via {GIT}}, {Osaka J. Math.}, {\bf 45} (2008), no. 4, {1027--1041}.

\bibitem{schweizer} A. Schweizer, {\em Some remarks on bielliptic and trigonal curves}, arXiv:1512.07963 [math.AG]


\bibitem{gon1} N.~Wangyu, {\em Cyclic covering of the projective line with prime gonality}, J. Pure and Applied Algebra {\bf 219} (2015), 1704--1710. 

\bibitem{gon2} N.~Wangyu, F.~Sakai, {\em Hyperelliptic curves among cyclic coverings of the projective line, II}, Arch. Math. {\bf 102} (2) (2014) 113--116.

\bibitem{Xiao} G.~Xiao, {\em Fibered algebraic surface with low slope}, Math. Ann. {\bf 276}, (1987), 449--466.

\bibitem{xiao-irregular} G.~Xiao, {\em Irregularity of surfaces with a linear pencil}, {Duke Math. J.}, {\bf 55} (1987), no. {3}, {597--602}.



\end{thebibliography}
\end{document}